\documentclass[a4paper,12pt]{article}
\usepackage{amsfonts,amsmath,amssymb,amsthm}
\usepackage{txfonts}
\usepackage[utf8]{inputenc}
\usepackage{amsmath}
\usepackage{amssymb}
\usepackage{mathtools}
\usepackage{bbm}

\RequirePackage[colorlinks,citecolor=red,urlcolor=blue, linkcolor=blue]{hyperref}
\usepackage[margin=0.8in]{geometry}

\usepackage[dvips]{graphics}
\usepackage{epsfig,rotating}
\usepackage{tabularx}
\usepackage{bm}
\usepackage{bbm}
\usepackage{float}
\usepackage{mdframed}
\usepackage{lipsum}
\usepackage{enumitem}
\usepackage{tikz}

\newtheorem{remark}{Remark}%

\usepackage[nottoc,numbib]{tocbibind} 
\usepackage[toc,page]{appendix}

\newtheorem{thm}{Theorem}

\newcommand{\beaa}{\begin{eqnarray*}}
\newcommand{\eeaa}{\end{eqnarray*}}
\newcommand{\bea}{\begin{eqnarray}}
\newcommand{\eea}{\end{eqnarray}}

\newcommand{\la}{\left\{}
\newcommand{\ra}{\right\}}
\newcommand{\lb}{\left(}
\newcommand{\rb}{\right)}
\newcommand{\mb}{\mathbb}

\begin{document}
\title{A note on time-uniform concentration inequality for matrix products}
\author{Tuan Pham \\ \small{Department of Statistics and Data Science, University of Texas, Austin} \\  \small{\href{mailto:tuan.pham@utexas.edu}{tuan.pham@utexas.edu} }
\and Alessandro Rinaldo \\ \small{Department of Statistics and Data Science, University of Texas, Austin} \\  \small{\href{mailto:alessandro.rinaldo@austin.utexas.edu}{alessandro.rinaldo@austin.utexas.edu}}
}
\date{}

\maketitle

\begin{abstract}

This short note contains a simple argument that allows us to go from fixed-time to any-time bound for the concentration of matrix products. The result presented here is motivated by the analysis of the Oja's algorithms. 

\end{abstract}

\section{Introduction}

Let $\bm{X}_1, \bm{X}_2,\dots, \bm{X}_n$ be a sequence of  i.i.d. positive-semidefinite (PSD) random matrices with $\mb{E} \bm{X}_1 =  \bm{\Sigma}$. From now on, unless stated otherwise, $\|.\|$ denotes the operator norm of matrices.   The goal of this short note is to prove a time-uniform concentration inequality for the matrix product
\[
\bm{Z}_n : = \prod_{i=n}^{1} \lb \bm{I}_d + \eta_{i}\bm{X}_i \rb.  
\]
In other words, we want a bound of the form 
\begin{align} \label{time uniform}
\mb{P} \lb \| \bm{Z}_n - \mb{E}\bm{Z}_n   \| \geq t\lb \delta,n \rb  \rb \leq \delta 
\end{align}
for some threshold $t(\delta,n)$ depending on $n$ and $\delta$.

The fixed-time concentration analog of \eqref{time uniform} has been studied extensively in the literature. The first result establishing a concentration inequality for $\bm{Z}_n$ appears in \cite{henriksen2020concentration,kathuria2020concentration}, where the authors derive such an inequality under the assumption that $\| \bm{X}_i - \bm{\Sigma} \|$ is bounded. Their argument is based on expanding the matrix product and partitioning it into blocks of independent terms, followed by a union bound. Subsequent improvements were developed in \cite{huang2022matrix} (see also the references therein), where the authors exploit the uniform smoothness of the Schatten $p$-norms.

The goal of this note is to provide an any-time bound for $\bm{Z}_n$ using a simple argument based on submartingale techniques, albeit at the cost of an additional sub-optimal term of order $\mathbb{E}\bm{Z}_n$.

\section{Settings and results}

Define 
\begin{align*}
    \bm{E}_n :&= \mb{E} \bm{Z}_n; \\
    m_i:&= \| \bm{I}_d + \eta_i \bm{\Sigma} \|; \\
    M_k:&= \prod_{i=1}^k m_i.
\end{align*}
and 
\begin{align*}
    \sigma_i:&= \eta_i \cdot \| \bm{X}_i - \bm{\Sigma} \|; \\
    V_k :&= \sum_{i=1}^k \sigma_i^2. 
\end{align*}
Suppose the step sizes $\eta_i$'s are chosen such that
\begin{align} \label{<= e}
M_n \sqrt{2 V_n \log(d/\delta)} \leq 1.
\end{align}
for some $\delta>0$. Our main result can be stated as follows.
\begin{thm} \label{any time}
    For all $\delta>0, \eta >1$ and all functions $h: \mb{N} \to \mb{R}_+$ such that 
    \[
    \sum_{k=0} \frac{1}{h(k)} \leq 1,
    \]
    we have 
    \[
    \mb{P} \lb  \exists n \colon \| \bm{Z}_n - \bm{E}_n \| \geq t\lb \frac{\delta}{h(k_n)}, \lfloor \eta^{k_n+1} \rfloor \rb  \rb \leq \delta
    \]
   whenever \eqref{<= e} holds, where
    \begin{align}
        k_n:&= \min \la  k=0,1,\ldots, \colon \lceil \eta^k \rceil \leq n \leq \lfloor \eta^{k+1} \rfloor \ra;   \label{k} \\
        t(\delta,n):&= e \| \bm{E}_n \| M_n \sqrt{2 V_n \log(d/\delta)}. \label{t}
    \end{align}
\end{thm}

\noindent \textbf{Proof of Theorem \ref{any time}.} Unlike the existing results that exploit supermartingale arguments, we will use a submartingale argument instead. The argument presented below is very simple, but is sub-optimal by a factor of order $\| \bm{E}_n \|$. Let us split the proof into a few steps. 

\textit{ \underline{Step 1: Constructing the submartingale.}}  Observe that 
\[
\bm{Y}_n = \bm{E}_n^{-1} \bm{Z}_n - \bm{I}_d 
\]
is a martingale with repsect to the natural filtration. 

By using the convexity of the map $\bm{X} \to \| \bm{X} \|$, we deduce that the process 
\[
\la \| \bm{Y}_k \|; k \geq 1 \ra
\]
is a submartingale. 

Moreover, we have the two-sided bound
\[
\| \bm{Y}_k \| \leq \| \bm{Z}_k - \bm{E}_k \| \leq \| \bm{E}_k \| \times \| \bm{Y}_k \|.
\]
for all $k \geq 1$.

The first inequality in the above is true since
\[
\| \bm{Y}_k \| \leq \| \bm{E}_k^{-1} \| \times \| \bm{Z}_k - \bm{E}_k \| \leq \underbrace{\prod_{i=1}^k \| \lb \bm{I}_d + \eta_i \bm{\Sigma} \rb^{-1} \|}_{\leq 1} \times \| \bm{Z}_k - \bm{E}_k \|.
\]
The second inequality follows from the simple observation that
\[
\| \bm{Z}_k - \bm{E}_k \| = \|  \bm{E}_k \bm{Y}_k  \| \leq \| \bm{E}_k \| \times \| \bm{Y}_k \|.
\]
Thus,
\begin{align} \label{any time 2}
    \mb{P} \lb \forall n \geq 1: \| \bm{Z}_n - \bm{E}_k \| \geq r_n \rb \leq \mb{P} \lb \forall n \geq 1: \| \bm{Y}_n \| \geq \frac{r_n}{\| \bm{E}_n \|} \rb.
\end{align}

\textit{ \underline{Step 2: Bounding the submartingale.}} To bound the submartingale in the RHS of \eqref{any time 2}, let us first divide $\mb{N}$ into intervals of the form
\[
\Big[\lfloor \eta^{k} \rfloor, \lfloor \eta^{k+1} \rfloor \Big]
\]
for $\eta>1$ as specified in the statement of Theorem \ref{any time}.

With $k_n$ and $t(\delta,n)$ as in \eqref{k} and \eqref{t}, repsectively, we have 
\begin{align*}
 \mathbb{P} \lb \exists n \colon \| \bm{Y}_n \| \geq t \lb \frac{\delta}{h(k_n)}, \lfloor \eta^{k_n+1} \rfloor \rb \rb
    & \leq  \mathbb{P} \lb \exists k, n \colon \lceil \eta^k \rceil \leq n \leq \lfloor \eta^{k+1} \rfloor,  \| \bm{Y}_n \| \geq t \lb \frac{\delta}{h(k_n)}, \lfloor \eta^{k_n+1} \rfloor  \rb \rb\\
    & = \mathbb{P} \lb \exists k, n \colon \lceil \eta^k \rceil \leq n \leq \lfloor \eta^{k+1} \rfloor,  \| \bm{Y}_n \| \geq t\lb \frac{\delta}{h(k)}, \lfloor \eta^{k+1} \rfloor \rb \rb \\
    & \leq \sum_{k=0}^\infty \mathbb{P} \lb \max_{ n \leq \lfloor \eta^{k+1} \rfloor} \| \bm{Y}_n \| \geq t\lb \frac{\delta}{h(k)}, \lfloor \eta^{k+1} \rfloor \rb   \rb.
    \end{align*}
For any $n\geq 1$ and $t>0$, by using the Doob's maximal inequality, we have 
\begin{align*}
    \mb{P} \lb \max_{k \leq n} \| \bm{Y}_n \| \geq u \rb \leq \min_{p \geq 1} \frac{\mathbb{E}[\max_{k \leq n}  \| \bm{Y}_k \|^p ]}{u^p} \leq \min_{p \geq 1} \frac{\mathbb{E}[  \| \bm{Y}_n \|^p ]}{u^p} \leq  \min_{p \geq 1} \frac{\mathbb{E}[  \| \bm{Z}_n - \bm{E}_n \|^p]}{u^p}.
\end{align*}
To bound the last term, one can follow the proof of Corollary 5.6 in \cite{huang2022matrix} to deduce that
\[
\mathbb{P} \left( \max_{k \leq n} \| \bm{Y}_k \| \geq u M_n \right)  \leq \max\{d,e\} \times \exp \left( -\frac{u^2}{2 e^2 V_n} \right), \quad \text{for} \ u \in [0,e].
\]
Thus, by choosing $t$ as in \eqref{t}, we have
\[
 \mathbb{P} \lb \max_{ n \leq \lfloor \eta^{k+1} \rfloor} \| \bm{Y}_n \| \geq t\lb \frac{\delta}{h(k)}, \lfloor \eta^{k+1} \rfloor \rb   \rb \leq \frac{\delta}{h(k)}.
\]
The proof is completed by employing the fact that $\sum_{k \geq 1} 1/h(k) \leq 1.$

\begin{remark}
    The boundary chosen in Theorem \ref{any time} is piece-wise constant over geometrically increasing epochs. One can easily turn it into a smooth boundary by choosing 
    \[
    f(n) =  e \| \bm{E}_n \|_{\mathrm{op}} M_n \sqrt{V_n (\log(d) + \log(\zeta(\alpha)/\delta) + \alpha \log( \log_\eta (n)  + 1)}
    \]
    where $\min \la \eta, \alpha \ra>1$ and $h(k) = (k+1)^\alpha \zeta(\alpha)$, with $\zeta(\cdot)$ the Riemann zeta function.

    Under the assumption \eqref{<= e}, the thresholds $\la f(n); n \geq 1 \ra$ satisfies
    \[
    \mb{P} \lb \exists n \geq 1: \| \bm{Z}_n - \bm{E}_n \| \geq f(n) \rb \leq \delta.
    \]
\end{remark}

\begin{remark}
    Regarding the tightness Theorem \ref{any time}, compared to the point-wise (in $n$) concentration bound of Huang et al. (2021), which is of order $e  M_n \sqrt{2 V_n \log(d/\delta)} $, we have picked up an additional term, namely $\| E_n \|_{\mathrm{op}}$. For our problem, this can be replaced by $M_n$. 
    
    To exemplify, let's consider the bound at  fixed time $n$. Set $\lambda_{\max}(\Sigma) = \mu$ and assume bounded data within an Euclidean ball of radius $L$:
    \[
    \| \bm{X}_i - \bm{\Sigma} \| \leq L.
    \]
    Also, set $\eta_i = 1/n$ for all $i \leq n$. Then, the bound of Huang et al. (2021) is of order 
\[
L e^\mu\sqrt{ \frac{ \log(d/\delta)}{n} }
\]
while our uniform bound, in addition to the extra iterated $\log$ terms, will be of order
\[
L e^{2\mu}\sqrt{ \frac{ \log(d/\delta)}{n} }.
\]
\end{remark}

\bibliographystyle{plain}
\bibliography{references}

\end{document}